\documentclass{birkjour}
\usepackage{amsthm}
\usepackage{mathtools}
\usepackage{url}
\usepackage{helvet}
\usepackage{times}
\usepackage{amsmath}
\usepackage{amsfonts}
\usepackage{amssymb}
\usepackage{mathrsfs}
\usepackage{bbm}
\usepackage{floatflt}
\usepackage{graphics}
 \newtheorem{thm}{Theorem}[section]
 
 \newtheorem{lem}[thm]{Lemma}
 
 \theoremstyle{definition}
 \newtheorem{defn}[thm]{Definition}
 \theoremstyle{remark}

 \numberwithin{equation}{section}
\begin{document}
\title{Super ${\bf 3}$-Lie Algebras Induced by Super Lie Algebras}
\author{Viktor Abramov}
\address{Institute of Mathematics, University of Tartu\\
Liivi 2 -- 602, Tartu 50409, Estonia}
\email{viktor.abramov@ut.ee}


\date{}

\begin{abstract}
We propose a notion of a super $n$-Lie algebra and construct a super $n$-Lie algebra with the help of a given binary super Lie algebra which is equipped with an analog of a supertrace. We apply this approach to the super Lie algebra of a Clifford algebra with even number of generators and making use of a matrix representation of this super Lie algebra given by a supermodule of spinors we construct a series of super 3-Lie algebras labeled by positive even integers.
\end{abstract}

\maketitle
\section{Introduction}
In search for new fundamental structures in theoretical physics and mechanics the physicists are increasingly turning their attention to algebraic structures which are based on ternary multiplication law (more generally on $n$-ary multiplication law). In 1973 Y. Nambu proposed a generalization of Hamiltonian mechanics where he replaced a canonically conjugate pair of variables by a triple of canonical variables and the usual Poisson bracket by a ternary operation (Nambu bracket) \cite{Nambu}. A geometric formalism for this generalization of Hamiltonian mechanics based on a notion of Nambu bracket of order $n$, the fundamental identity and a concept of Nambu-Poisson manifold was developed by L. Takhtajan in \cite{Takhtajan}. An important part of this geometric formalism is a notion of an $n$-Lie algebra which was also studied by V. T. Filippov in \cite{Filippov}. J. Arnlind, A. Makhlouf and S. Silvestrov constructed an $n$-Lie algebra by means of a binary Lie algebra endowed with an analog of a trace \cite{Makhlouf-Silvestrov 2} and this $n$-Lie algebra was called an $n$-Lie algebra induced by Lie algebra. They also studied the cohomologies of $n$-Lie algebra induced by a Lie algebra and found the relation of these cohomologies to the cohomology of initial binary Lie algebra. In this paper we propose a notion of super $n$-Lie algebra and follow an approach of J. Arnlind, A. Makhlouf and S. Silvestrov to construct a super $n$-Lie algebra with the help of a binary super Lie algebra which is equipped with an analog of a supertrace, and this super $n$-Lie algebra is called the super $n$-Lie algebra induced by a super Lie algebra. We apply this approach to the super Lie algebra of Clifford algebra (with even number of generators) with the matrix representation given by a supermodule of spinors and for any even integer $n\geq 2$ we construct a super 3-Lie algebra.
\section{$\bf n$-Lie algebras}
In this section we remind the definition of $n$-Lie algebra.

Let ${\mathbb K}$ be a field of real or complex numbers, $\frak g$ be a vector space over $\mathbb K$, and ${\frak g}^n={\frak g}\times{\frak g}\times\ldots\times{\frak g} (n \mbox{ times})$. A vector space $\mathfrak g$ is said to be an $n$-Lie algebra if $\mathfrak g$ is endowed with a multilinear skew-symmetric mapping
$[\cdot,\cdot,\ldots,\cdot]:{\mathfrak g}^n\to {\mathfrak g}$ which satisfies the identity
\begin{eqnarray}
&&[y_1,y_2,\ldots,y_{n-1},[x_1,x_2,\ldots,x_n]] = \nonumber\\
  &&\qquad\quad\sum_{k=1}^{n}[x_1,x_2,\ldots,x_{k-1},[y_1,y_2,\ldots,y_{n-1},x_{k}],x_{k+1},\ldots,x_n],
  \label{identity for n-Lie algebra}
\end{eqnarray}
where $y_1,y_2,\ldots,y_{n-1},x_1,x_2,\ldots,x_n\in{\frak g}$. Particularly if $n=2$ then the above definition yields the definition of a (binary) Lie algebra and the identity (\ref{identity for n-Lie algebra}) takes on the form of Jacobi identity. Let $\frak h$ be a (binary) Lie algebra with a Lie bracket $[\cdot,\cdot]:{\frak h}\times{\frak h}\to {\frak h}$, and $\phi:{\frak h}^k\to {\mathbb K}, k\geq 1$ be an $\mathbb K$-multilinear skew-symmetric $k$-form which satisfies the condition
\begin{equation}
\sum_{m=1}^k \phi(x_1,x_2,\ldots,x_{m-1},[x_m,y],x_{m+1},\ldots,x_k)=0,
\end{equation}
where $x_1,x_2,\ldots,x_k,y$ are elements of $\frak h$. Particularly if $k=1$ then a form $\phi:{\frak h}\to {\mathbb K}$ satisfies $\phi([x,y])=0$ for any $x,y\in{\frak h}$.
\begin{thm}
Let $k\geq 1$ be an integer, $\sigma=(i_1,i_2,\ldots,i_k,i_{k+1},i_{k+2})$ be a permutation of the integers $(1,2,\ldots,k+2)$ such that $1\leq i_1<i_2,\ldots<i_k\leq k+2$, $1\leq i_{k+1}<i_{k+2}\leq k+2$ and $|\sigma|$ be the parity of $\sigma$. Define
\begin{equation}
[x_1,x_2,\ldots,x_{k+2}]=\sum_{\sigma} (-1)^{|\sigma|}\phi(x_{i_1},x_{i_2},\ldots,x_{i_k})[x_{i_{k+1}},x_{i_{k+2}}].
\label{n-bracket for Lie algebra}
\end{equation}
Then $\frak h$ endowed with (\ref{n-bracket for Lie algebra}) is the $n$-Lie algebra, where $n=k+2$.
\end{thm}
It is worth noting that in the particular case of $k=1$ this theorem yields the 3-Lie algebra with ternary Lie bracket
\begin{equation}
[x,y,z]=\phi(x)[y,z]+\phi(y)[z,x]+\phi(z)[x,y],
\label{ternary Lie bracket}
\end{equation}
and a form $\phi:{\frak h}\to {\mathbb K}$ which for any $x,y\in{\frak h}$ satisfies $\phi([x,y])=0$ can be viewed as an analog of a trace. The 3-Lie algebras of this kind induced by a binary Lie algebra are introduced and studied \cite{Makhlouf-Silvestrov}.

\section{Super $\bf n$-Lie algebras}
In this section we give a definition of a super $n$-Lie algebra and prove that one can construct a super $n$-Lie algebra by means of a super trace.

Let $V=V_0\oplus V_1$ be a finite-dimensional super vector space, where $V_0$ is the subspace of even elements and $V_1$ is the subspace of odd elements. If $v\in V$ is a homogeneous element then its degree will be denoted by $|v|$, where $|v|\in {\mathbb Z}_2$ and ${\mathbb Z}_2=\{0,1\}$. Let $V^n=V\times V\times\ldots\times V$. For any ${\bf v}=(v_1,v_2,\dots,v_n)\in V^n$, where $v_1,v_2,\dots,v_n$ are homogeneous elements, and for every integer $k=1,2,\ldots,n$ we define
$$
|{\bf v}|_k=\sum_{i=1}^k|v_i|.
$$
Let $\mbox{End}\,(V)$ be the super vector space of endomorphisms of a super vector space $V$. The composition of two endomorphisms $a\circ b$  determines the structure of superalgebra in $\mbox{End}\,(V)$, and the graded binary commutator $[a,b]=a\circ b-(-1)^{|a||b|}$ induces the structure of super Lie algebra in $\mbox{End}\,(V)$. The supertrace of an endomorphism $a:V\to V$ can be defined by
\begin{equation}
\mbox{Str}(a)=\left\{
\begin{array}{ll}
\mbox{Tr}(a|_{V_0})-\mbox{Tr}(a|_{V_1}) & \text{if } a \text{ is even},\\
0 & \text{if } a \text{ is odd} .
\end{array} \right.
\end{equation}
For any endomorphisms $a,b$ it holds $\mbox{Str}([a,b])=0$.

Given ${\bf v}=(v_1,v_2,\dots,v_n)\in V^n$ and a permutation $\sigma=(i_1,i_2,\ldots,i_n)$ of the integers $(1,2,\ldots,n)$ we assign to them the element ${\bf v}_{\sigma}\in V^n$ defined by
$$
({\bf v}, {\sigma})\;\rightarrow\;{\bf v}_{\sigma}=(v_{i_1},v_{i_2},\ldots,v_{i_n}).
$$
Let $i_k$ be an element of permutation $\sigma$. Assume that $(i_{k_1},i_{k_2},\ldots,i_{k_r})$ are the elements of permutation $\sigma$ which precede an element $i_k$ in $\sigma$ and greater than $i_k$. Obviously $r$ is the number of inversions of an element $i_k$ in $\sigma$. Define the integer $|{\bf v}_\sigma|\in{\mathbb Z}_2$ by
$$
|{\bf v}_\sigma|=\sum_{k=1}^n |v_{i_k}|(|v_{i_{k_1}}|+|v_{i_{k_2}}|+\ldots+|v_{i_{k_r}}|).
$$
\begin{defn}
A super vector space $\frak g={\frak g}_0\oplus {\frak g}_1$ is said to be a super $n$-Lie algebra if it is endowed with an $n$-ary bracket $[ \cdot , \cdot ,\ldots,\cdot]: {\frak g}\times {\frak g}\times\ldots\times{\frak g}($n$\;\mbox{times})\to {\frak g}$ which satisfies the following conditions:
\begin{enumerate}
\item
 $n$-ary bracket is a multilinear mapping and the degree of $n$-ary bracket of $n$ homogeneous elements is equal to the sum of degrees of these elements, i.e. for any homogeneous element ${\bf x}=(x_1,x_2,\ldots,x_n)\in {\frak g}^n$ it holds
 $$
 |[x_1,x_2,\ldots,x_n]|=|{\bf x}|_n,
 $$
\item
 $n$-ary bracket is a graded skew-symmetric multilinear mapping, i.e. for any homogeneous elements $x_1,x_2,\ldots,x_n\in {\frak g}$ and for any $k=1,2,\ldots,n-1$ it holds
$$
[x_1,x_2,\ldots,x_k,x_{k+1},\ldots,x_n]=-(-1)^{|x_k||x_{k+1}|}[x_1,x_2,\ldots,x_{k+1},x_{k},\ldots,x_n],
$$
\item
 $n$-ary bracket satisfies the identity
\begin{eqnarray}
&&[y_1,y_2,\ldots,y_{n-1},[x_1,x_2,\ldots,x_n]] = \nonumber\\
  &&\;\;\sum_{k=1}^{n}(-1)^{|{\bf x}|_{k-1}|{\bf y}|_{n-1}}[x_1,x_2,\ldots,x_{k-1},[y_1,y_2,\ldots,y_{n-1},x_{k}],x_{k+1},\ldots,x_n],\nonumber
\end{eqnarray}
 where ${\bf y}=(y_1,y_2,\ldots,y_{n-1})\in {\frak g}^{n-1}, {\bf x}=(x_1,x_2,\ldots,x_n)\in{\frak g}^n$ are homogeneous elements.
\end{enumerate}
\end{defn}
\begin{lem}
Let $\mbox{End}\,(V)$ be the super vector space of endomorphisms of a super vector space $V$ and ${\bf a}=(a_1,a_2,\ldots,a_n)$ be a sequence of $n$ endomorphisms of $V$. Define
\begin{equation}
[a_1,a_2,\ldots,a_n]=\sum_{\sigma}(-1)^{|\sigma|+|\bf{a}_\sigma|}a_{i_1}\circ a_{i_2}\circ\ldots\circ a_{i_n},
\label{n-bracket for endomorphisms}
\end{equation}
where $\sigma=(i_1,i_2,\ldots,i_n)$ is a permutation of integers $(1,2,\ldots,n)$ and $|\sigma|$ is the parity of this permutation. Then the super vector space $\mbox{End}\,(V)$ endowed with the $n$-ary bracket (\ref{n-bracket for endomorphisms}) is the super $n$-Lie algebra.
\end{lem}
\begin{defn}
A representation of a super $n$-Lie algebra $\frak g$ is a linear mapping $\rho:{\frak g}\to \mbox{End}\,(V)$, where $V$ is a super vector space (a representation space of $\frak g$), which satisfies:
\begin{enumerate}
\item
for any homogeneous element $x\in{\frak g}$ its image $\rho(x)$ in $\mbox{End}\,(V)$ is the homogeneous element and $|\rho(x)|=|x|$,
\item
for any elements $x_1,x_2,\ldots,x_n\in {\frak g}$ it holds
$$
\rho([x_1,x_2,\ldots,x_n])=[\rho(x_1),\rho(x_2),\ldots,\rho(x_n)].
$$
\end{enumerate}
\end{defn}
\noindent
Let ${\frak g}$ be a super $n$-Lie algebra with $n$-ary bracket $[\cdot,\cdot,\ldots,\cdot]$, and $\rho:{\frak g}\to \mbox{End}\,(V)$ be a representation of $\frak g$. For any homogeneous element
${\bf x}=(x_1,x_2,\ldots,x_{n+1})$ of $\frak g$ we define the $(n+1)$-ary bracket by
\begin{eqnarray}
[x_1,x_2,\ldots,x_{n+1}]_\rho &=&\sum_{k=1}^{n+1}(-1)^{k-1}(-1)^{|x_k||{\bf x}|_{k-1}}\mbox{Str}(\rho(x_k))\nonumber\\
        &&\qquad\qquad\qquad\quad\times [x_1,x_2,\ldots,\hat x_k,\ldots,x_{n+1}],
\label{n-bracket with supertrace}
\end{eqnarray}
where the hat over an element $\hat x_k$ means that this element is omitted.
\begin{thm}
A super $n$-Lie algebra $\frak g$ equipped with the $(n+1)$-ary bracket (\ref{n-bracket with supertrace}) is the super $(n+1)$-ary Lie algebra.
\label{theorem for n+1 Lie algebra with super trace}
\end{thm}
\section{Supermodule over Clifford algebra}
A Clifford algebra is a unital associative algebra, which can be equipped with an $\mathbb Z_2$-graded structure, and it provides a well known example of a superalgebra. Taking the graded commutator of two elements of this superalgebra one can consider it as the super Lie algebra. In this section  we consider a supermodule over Clifford algebra with even number of generators, and this supermodule provides us with a representation of mentioned above super Lie algebra. Making use of a supertrace of this representation we construct a super 3-Lie algebra with the help of the formula (\ref{n-bracket with supertrace}).

We remind that Clifford algebra $C_n$ is the unital associative algebra over $\mathbb C$ generated by $\gamma_1,\gamma_2,\ldots,\gamma_n$ which obey the relations
\begin{equation}
\gamma_i\gamma_j+\gamma_j\gamma_i=2\,\delta_{ij}\,e,\quad i,j=1,2,\ldots,n,
\end{equation}
where $e$ is the unit element of Clifford algebra. Let ${\mathcal N}=\{1,2,\ldots,n\}$ be the set of integers from 1 to $n$. If $I$ is a subset of $\mathcal N$, i.e. $I=\{i_1,i_2,\ldots,i_k\}$ where $1\leq i_1<i_2<\ldots<i_k\leq n$, then one can associate to this subset $I$ the monomial $\gamma_I=\gamma_{i_1}\gamma_{i_2}\ldots\gamma_{i_k}$. If $I=\emptyset$ one defines $\gamma_{\emptyset}=e$. The number of elements of a subset $I$ will be denoted by $|I|$. It is obvious that the vector space of Clifford algebra $C_n$ is spanned by the monomials $\gamma_I$, where $I\subseteq {\mathcal N}$. Hence the dimension of this vector space is $2^n$ and any element $x\in C_n$ can be expressed in terms of these monomials as
$$
x=\sum_{I\subseteq {\mathcal N}} a_I\gamma_I,
$$
where $a_I=a_{i_1i_2\ldots i_k}$ is a complex number. It is easy to see that one can endow a Clifford algebra $C_n$ with the ${\mathbb Z}_2$-graded structure by assigning the degree $|\gamma_I|=|I|\,(\mbox{mod}\,2)$ to monomial $\gamma_I$. Then a Clifford algebra $C_n$ can be considered as the superalgebra since for any two monomials it holds $|\gamma_I\gamma_J|=|\gamma_I|+|\gamma_J|$.

Another way to construct this superalgebra which does not contain explicit reference to Clifford algebra is given by the following theorem.
\begin{thm}
Let $I$ be a subset of ${\mathcal N}=\{1,2,\ldots,n\}$, and $\gamma_I$ be a symbol associated to $I$. Let $C_n$ be the vector space spanned by the symbols $\gamma_I$. Define the degree of $\gamma_I$ by $|\gamma_I|=|I| (\,\mbox{mod}\,2)$, where $|I|$ is the number of elements of $I$, and the product of $\gamma_I,\gamma_J$ by
\begin{equation}
\gamma_I\,\gamma_J=(-1)^{\sigma(I,J)}\gamma_{I\Delta J},
\label{superalgebra product}
\end{equation}
where $\sigma(I,J)=\sum_{j\in J}\sigma(I,j)$, $\sigma(I,j)$ is the number of elements of $I$ which are greater than $j\in J$, and $I\Delta J$ is the symmetric difference of two subsets. Then $C_n$ is the unital associative superalgebra, where the unit element $e$ is $\gamma_\emptyset$.
\end{thm}
This theorem can be proved by means of the properties of symmetric difference of two subsets. We remind a reader that the symmetric difference is commutative $I\oplus J=J\oplus I$, associative $(I\Delta J)\Delta K=I\Delta (J\Delta K)$ and $I\Delta \emptyset=\emptyset\Delta I$. The latter shows that $\gamma_{\emptyset}$ is the unit element of this superalgebra. The symmetric difference also satisfies $|I\Delta J|=|I|+|J|\;(\mbox{mod}\,2)$. Hence $C_n$ is the superalgebra.

The superalgebra $C_n$ can be considered as the super Lie algebra if for any two homogeneous elements $x,y$ of this superalgebra one introduces the graded commutator $[x,y]=xy-(-1)^{|x||y|}yx$ and extends it by linearity to a whole superalgebra $C_n$. We will denote this super Lie algebra by ${\frak C}_n$. Then $\{\gamma_I\}_{I\subseteq {\mathcal N}}$ are the generators of this super Lie algebra ${\frak C}_n$, and its structure is entirely determined by the graded commutators of $\gamma_I$. Then for any two generators $\gamma_I,\gamma_J$ we have
\begin{equation}
[\gamma_I,\gamma_J]=f(I,J)\;\gamma_{I\Delta J},
\label{binary commutators}
\end{equation}
where $f(I,J)$ is the integer-valued function of two subsets of $\mathcal N$ defined by
$$
f(I,J)=(-1)^{\sigma(I,J)}\big(1-(-1)^{|I\cap J|}\big),
$$
It is easy to verify that the degree of graded commutator is consistent with the degrees of generators, i.e. $[\gamma_I,\gamma_J]=|\gamma_I|+|\gamma_J|.$ Indeed the function $\sigma(I,J)$ satisfies
$$
\sigma(I,J)=|I||J|-|I\cap J|-\sigma(I,J),
$$
and
\begin{eqnarray}
f(J,I) &=& (-1)^{\sigma(J,I)}\big(1-(-1)^{|I\cap J|}\big)\nonumber\\
    &=& (-1)^{|I||J|-|I\cap J|-\sigma(I,J)}\big(1-(-1)^{|I\cap J|}\big)\nonumber\\
    &=& (-1)^{|I||J|}(-1)^{\sigma(I,J)}\big((-1)^{|I\cap J|}-1\big)
    =-(-1)^{|I||J|}f(I,J). \nonumber
\end{eqnarray}
Hence $[\gamma_I,\gamma_J]=-(-1)^{|I||J|}[\gamma_J,\gamma_I]$ which shows that the relation (\ref{binary commutators}) is consistent with the symmetries of graded commutator. It is obvious that if the intersection of subsets $I,J$ contains an even number of elements then $f(I,J)=0$, and the graded commutator of $\gamma_I,\gamma_J$ is trivial. Particularly if at least one of two subsets $I,J$ is the empty set then $f(I,J)=0$. Thus any graded commutator (\ref{binary commutators}) containing $e$ is trivial.

As an example, consider the super Lie algebra $\frak C_2$. Its underlying vector space is 4-dimensional and $\frak C_2$ is generated by two even degree generators $e,\gamma_{12}$ and two odd degree generators $\gamma_1,\gamma_2.$ The non-trivial relations of this super Lie algebra are given by
\begin{eqnarray}
[\gamma_1,\gamma_1]=[\gamma_2,\gamma_2]=2\,e,\;
[\gamma_1,\gamma_{12}]=2\,\gamma_2,\; [\gamma_2,\gamma_{12}]=-2\,\gamma_1.
\end{eqnarray}

Now we assume that $n=2m, m\geq 1$ is an even integer. The super Lie algebra ${\frak C}_n$ has a matrix representation which can be described as follows. Fix $n=2$ and identify the generators $\gamma_1,\gamma_2$ with the Pauli matrices $\sigma_1,\sigma_2$, i.e.
\begin{equation}
\gamma_1=\left(
           \begin{array}{cc}
             0 & 1 \\
             1 & 0 \\
           \end{array}
         \right),\quad
\gamma_2=\left(
           \begin{array}{cc}
             0 & -i \\
             i & 0 \\
           \end{array}
         \right).
\label{Pauli matrices 1,2}
\end{equation}
Then $\gamma_{12}=\gamma_1\gamma_2=i\,\sigma_3$ where
$$
\sigma_3=\left(
           \begin{array}{cc}
             1 & 0 \\
             0 & -1 \\
           \end{array}
         \right).
$$
Let $S^2$ be the 2-dimensional complex super vector space $\mathbb C^2$ with the odd degree operators (\ref{Pauli matrices 1,2}), where the $\mathbb Z_2$-graded structure of $S^2$ is determined by $\sigma_3=i^{-1}\gamma_{12}$. Then $C_2\simeq \mbox{End}\,(S^2)$, and $S^2$ can be considered as a supermodule over the superalgebra $C_2$. Let $S^n=S^2\otimes S^2\otimes\ldots \otimes S^2 (m-\,\mbox{times})$. Then $S^n$ can be viewed as a supermodule over the $m$-fold tensor product of $C_2$, which can be identified with $C_n$ by identifying $\gamma_1,\gamma_2$ in the $j$th factor with $\gamma_{2j-1},\gamma_{2j}$ in $C_n$. This $C_n$-supermodule $S^n$ is called the supermodule of spinors \cite{Quillen-Mathai}. Hence we have the matrix representation for the Clifford algebra $C_n$, and this matrix representation or supermodule of spinors allows one to consider the supertrace, and it can be proved \cite{Quillen-Mathai} that
\begin{equation}
\mbox{Str}(\gamma_I)=\left\{
\begin{array}{ll}
0 & \text{if } I<{\mathcal N},\\
(2i)^m & \text{if } I={\mathcal N}.
\end{array} \right.
\label{super trace}
\end{equation}
Now we have the super Lie algebra $\frak C_n$ with the graded commutator defined in (\ref{binary commutators}) and its matrix representation based on the supermodule of spinors. Hence we can construct a super 3-Lie algebra by making use of graded ternary commutator (\ref{n-bracket with supertrace}). Applying the formula (\ref{n-bracket with supertrace}) we define the graded ternary commutator for any triple $\gamma_I,\gamma_J,\gamma_K$ of elements of basis for $\frak C_n$ by
\begin{eqnarray}
&& [\gamma_I,\gamma_J,\gamma_K] = \mbox{Str}(\gamma_I)\,[\gamma_J,\gamma_K]-(-1)^{|I||J|}\mbox{Str}(\gamma_J)\, [\gamma_I,\gamma_K]\nonumber\\
&&\qquad\qquad\qquad\qquad\qquad\qquad\quad + (-1)^{|K|(|I|+|J|)}\mbox{Str}(\gamma_K)\,[\gamma_I,\gamma_J],
\label{ternary graded commutator with super trace}
\end{eqnarray}
where the binary graded commutator at the right-hand side of this formula is defined by (\ref{binary commutators}). According to Theorem \ref{theorem for n+1 Lie algebra with super trace} the vector space spanned by $\gamma_I, I\subset {\mathcal N}$ and equipped with the ternary graded commutator (\ref{ternary graded commutator with super trace}) is the super 3-Lie algebra which will be denoted by ${\frak C}^{(3)}_n$. Making use of (\ref{binary commutators}) we can write the expression at the right-hand side of the above formula in the form
\begin{eqnarray}
&& [\gamma_I,\gamma_J,\gamma_K] =f(J,K) \mbox{Str}(\gamma_I)\,\gamma_{J\Delta K}-(-1)^{|I||J|}f(I,K)\mbox{Str}(\gamma_J)\, \gamma_{I\Delta K}\nonumber\\
&&\quad\qquad\qquad\qquad\qquad\qquad\qquad\quad + (-1)^{|K|(|I|+|J|)}f(I,J)\mbox{Str}(\gamma_K)\,\gamma_{I\Delta J}.\nonumber
\end{eqnarray}
From the formula for supertrace (\ref{super trace}) it follows immediately that the above graded ternary commutator is trivial if none of subsets $\gamma_i,\gamma_J,\gamma_K$ is equal to ${\mathcal N}$. Similarly this graded ternary commutator is also trivial if all three subsets $I,J,K$ are equal to $\mathcal N$, i.e. $I=J=K={\mathcal N}$, or two of them are equal to
 ${\mathcal N}$.
\begin{thm}
The graded ternary commutators of the generators $\gamma_I,I\subseteq {\mathcal N}$ of the super 3-Lie algebra ${\frak C}^{(3)}_n$ are given by
\begin{equation}
[\gamma_I,\gamma_J,\gamma_K]=\left\{
\begin{array}{ll}
(2i)^mf(I,J)\gamma_{I\Delta J} & \text{if } I\neq {\mathcal N}, J\neq {\mathcal N}, K={\mathcal N},\\
0& \text{in all other cases }.
\end{array} \right.
\label{proposition}
\end{equation}
\end{thm}
\subsection*{Acknowledgment}
The author gratefully acknowledges that this work was financially supported by the Estonian Science Foundation under the research grant ETF9328 and this work was also financially supported by institutional research funding IUT20-57
of the Estonian Ministry of Education and Research.

\end{document}